\documentclass[preprint,12pt]{elsarticle}

\usepackage{graphicx}
\usepackage[cp1251]{inputenc}
\usepackage{amssymb}
\usepackage{amsthm}
\usepackage{cmap}

\usepackage{amsmath}

\biboptions{sort&compress}

\newcommand{\sysn}{\left\{\begin{array}{rcl}}
\newcommand{\sysk}{\end{array}\right.}

\newtheorem{theorem}{Theorem}[section]
\newtheorem{lemma}[theorem]{Lemma}

\theoremstyle{example}

\theoremstyle{definition}
\newtheorem{definition}[theorem]{Definition}

\journal{...}

\begin{document}

\title{Baireness of the space of  pointwise stabilizing functions of the first Baire class}

\author{Alexander V. Osipov}

\address{Krasovskii Institute of Mathematics and Mechanics, \\ Ural Federal
 University, Ural State University of Economics, Yekaterinburg, Russia}

\ead{OAB@list.ru}

\begin{abstract} A topological space $X$ is {\it Baire} if the Baire
Category Theorem holds for $X$, i.e., the intersection of any
sequence of open dense subsets of $X$ is dense in $X$. In this
paper, we have obtained that the space $B^{st}_1(X)$ of pointwise
stabilizing Baire-one functions is Baire if the space $B_1(X)$ of
Baire-one functions is so.
 This answers a question posed
recently by T. Banakh and S. Gabriyelyan.

\end{abstract}


\maketitle 


\section{Introduction}
For a topological space $X$, we denote by $C_p(X)$ the family
$C(X)$ of all continuous function from $X$ to the real line
$\mathbb{R}$ endowed with the topology of pointwise convergence.

A real-valued function $f$ on a space $X$ is a {\it Baire-one
function} (or a {\it function of the first Baire class}) if $f$ is
a pointwise limit of a sequence of continuous functions on $X$.
Let $B_1(X)$ denote the family of all Baire-one real-valued
functions on a space $X$ endowed with the topology of pointwise
convergence.

We say that a sequence $\{f_n\}_{n\in\omega}\subseteq Y^X$ {\it
stable (or discretely) converges} to a function $f\in Y^X$ if for
every $x\in X$ the set $\{n\in \omega: f_n(x)\neq f(x)\}$ is
finite.
 Let $B^{st}_1(X)$ be the
family of all functions $f:X\rightarrow \mathbb{R}$ which are
limits of sequence $\{f_n\}_{n\in\omega}\subseteq C_p(X)$ stably
converging to $f$. Functions in the family $B_1^{st}(X)$ are
called the {\it functions of stable first Baire class (or
pointwise stabilizing functions of the first Baire class)}.

 A space is {\it meager} (or {\it of the first Baire category}) if it
can be written as a countable union of closed sets with empty
interior. A topological space $X$ is {\it Baire} if the Baire
Category Theorem holds for $X$, i.e., the intersection of any
sequence of open dense subsets of $X$ is dense in $X$. Clearly, if
$X$ is Baire, then $X$ is not meager.

 Being a Baire space is an important topological
property for a space and it is therefore natural to ask when
function spaces are Baire. The Baire property for continuous
mappings was first considered in \protect\cite{Vid}. Then a paper
\protect\cite{ZMc} appeared, where various aspects of this topic
were considered. In \protect\cite{ZMc}, necessary and, in some
cases, sufficient conditions on a space $X$ were obtained under
which the space $C_p(X)$ is Baire. The problem about a
characterization of Baireness for $C_p(X)$ was solved
independently by Pytkeev \protect\cite{pyt1}, Tkachuk
\protect\cite{tk} and van Douwen \protect\cite{vD}.

\medskip

\begin{theorem}(Pytkeev-Tkachuk-van Douwen)  The space
$C_p(X)$ is Baire if, and only if, every pairwise disjoint
sequence of non-empty finite subsets of $X$ has a strongly
discrete subsequence.
\end{theorem}

A collection $\mathcal{G}$ of subsets of $X$ is {\it discrete} if
each point of $X$ has a neighborhood meeting at most one element
of $\mathcal{G}$, and is {\it strongly discrete} if for each $G\in
\mathcal{G}$ there is an open superset $U_G$ of $G$ such that
$\{U_G: G\in\mathcal{G}\}$ is discrete.

\medskip

A $Coz_{\delta}$-subset of $X$ containing $x$ is called a {\it
$Coz_{\delta}$ neighborhood} of $x$.

\begin{definition}\cite{Os}  A set $A\subseteq X$ is called {\it strongly
$Coz_{\delta}$-disjoint}, if there is a pairwise disjoint
collection $\{F_a: F_a$ is a $Coz_{\delta}$ neighborhood of $a$,
$a\in A\}$  such that $\{F_a: a\in A\}$ is {\it completely
$Coz_{\delta}$-additive}, i.e. $\bigcup\limits_{b\in B} F_b\in
Coz_{\delta}$ for each $B\subseteq A$.

A disjoint sequence $\{\Delta_n: n\in \mathbb{N}\}$ of (finite)
sets is  {\it strongly $Coz_{\delta}$-disjoint} if the set
$\bigcup\{\Delta_n: n\in \mathbb{N}\}$ is strongly
$Coz_{\delta}$-disjoint.

\end{definition}

In \cite{Os}, we have obtained a characterization when $B_1(X)$ is
Baire.

\begin{theorem}\label{th22} Let $X$ be a topological space. The following assertions are
equivalent:

1. $B_1(X)$ is Baire;

2. every pairwise disjoint sequence of non-empty finite subsets of
$X$ has a strongly $Coz_{\delta}$-disjoint subsequence.
\end{theorem}

\medskip

In \cite{BG}, T.Banakh and S.Gabriyelyan considers the following
question (Problem 6.10 and Problem 9.11): {\it Is there a
topological space $X$ such that $B_1(X)$ is Baire and
$B^{st}_1(X)$ is meager?}

\medskip

In this paper, we get the following theorem which answers the
Banakh-Gabriyelyan question.

\begin{theorem}\label{th23} Let $X$ be a topological space. The following assertions are
equivalent:

1. $B_1(X)$ is Baire;

2. $B^{st}_1(X)$ is Baire;

3. Every pairwise disjoint sequence $\{\Delta_n: n\in
\mathbb{N}\}$ of non-empty finite subsets of $X$ has a subsequence
$\{\Delta_{n_k}: k\in \mathbb{N}\}$  such that there is a pairwise
disjoint family $\{F_a: F_a$ is a zero-set neighborhood of $a$,
$a\in \bigcup \Delta_{n_k}\}$  and $\{F_a: a\in \bigcup
\Delta_{n_k}\}$ is completely $Coz_{\delta}$-additive.
\end{theorem}

\section{Main definitions and notation}

Throughout this paper,  all spaces are assumed to be Tychonoff.
 The set of positive integers is denoted by $\mathbb{N}$ and
$\omega=\mathbb{N}\cup \{0\}$. Let $\mathbb{R}$ be the real line,
we put $\mathbb{I}=[0,1]\subset \mathbb{R}$, and let $\mathbb{Q}$
be the rational numbers. Let $f:X\rightarrow \mathbb{R}$ be a
real-valued function, then $\parallel f \parallel= \sup \{|f(x)|:
x\in X\}$, $S(g,\epsilon)=\{f: \parallel g-f
\parallel<\epsilon\}$, $B(g,\epsilon)=\{f: \parallel g-f
\parallel\leq\epsilon\}$, where $g$ is a real valued function and
$\epsilon>0$. Let $V=\{f\in \mathbb{R}^X: f(x_i)\in V_i,
i=1,...,n\}$ where $x_i\in X$, $V_i\subseteq \mathbb{R}$ are
bounded intervals for $i=1,...,n$, then $supp V=\{x_1,...,x_n\}$ ,
$diam V=\max \{diam V_i : 1\leq i \leq n \}$.

 We recall that a subset of $X$ that is the
 complete preimage of zero for a certain function from~$C(X)$ is called a {\it zero-set}.
A subset $O\subseteq X$  is called  a cozero-set (or functionally
open) of $X$ if $X\setminus O$ is a zero-set of $X$. It is easy to
check that zero sets are preserved by finite unions and countable
intersections. Hence cozero sets are preserved by finite
intersections and countable unions. Countable unions of zero sets
will be denoted by $Zer_{\sigma}$ (or $Zer_{\sigma}(X)$),
countable intersection of cozero sets by $Coz_{\delta}$ (or
$Coz_{\delta}(X)$). It is easy to check that $Zer_{\sigma}$-sets
are preserved by countable unions and finite intersections.  Note
that any zero set is $Coz_{\delta}$ and any cozero-set is
$Zer_{\sigma}$. It is well known that $f$ is of the first Baire
class if and only if $f^{-1}(U)\in Zer_{\sigma}$ for every open
$U\subseteq \mathbb{R}$ (see Exercise 3.A.1 in
\protect\cite{lmz}).

\medskip

\begin{lemma}\label{lem101}(Lemma 3.3 in \cite{Os}) Let $X$ be a Hausdorff space  and let
$\{F_i: i\in \mathbb{N}\}$ forms a disjoint completely
$Coz_{\delta}$-additive system. Then any family $\{L_i:
L_i\subseteq F_i$, $i\in\mathbb{N}\}$ of $Coz_{\delta}$ subsets of
$X$ is a completely $Coz_{\delta}$-additive system.
\end{lemma}

\section{Proof of Theorem 1.4}

\begin{proof}
$(2)\Rightarrow(1)$. Since $B^{st}_1(X)$ is a dense subset of
$B_1(X)$, then $B_1(X)$ is Baire if $B^{st}_1(X)$ is so.

\medskip

$(1)\Rightarrow(3)$. By Theorem \ref{th22}, every pairwise
disjoint sequence of non-empty finite subsets of $X$ has a
strongly $Coz_{\delta}$-disjoint subsequence. Let $\{\Delta_n:
n\in \mathbb{N}\}$ be a pairwise disjoint sequence of non-empty
finite subsets of $X$ and $\{\Delta_{n_k}: k\in \mathbb{N}\}$ be a
strongly $Coz_{\delta}$-disjoint subsequence. Then, there is a
pairwise disjoint collection $\{F_a: F_a$ is a $Coz_{\delta}$
neighborhood of $a$, $a\in  \bigcup_{k\in \mathbb{N}}
\Delta_{n_k}\}$  such that $\{F_a: a\in \bigcup_{k\in \mathbb{N}}
\Delta_{n_k}\}$ is {\it completely $Coz_{\delta}$-additive}. For
every $a\in F_a$ there is a zero-set $W_a$ of $X$ such that $a\in
W_a\subseteq F_a$. By Lemma \ref{lem101}, the family $\{W_a: a\in
\bigcup_{k\in \mathbb{N}} \Delta_{n_k}\}$ is completely
$Coz_{\delta}$-additive.

\medskip

$(3)\Rightarrow(2)$. Assume contrary. Let
$B^{st}_1(X)=\bigcup_{n\in\mathbb{N}} F_n$, where $F_n$ is nowhere
dense in $B_1(X)$ and $F_n\subseteq F_{n+1}$ for every $n\in
\mathbb{N}$.
\medskip

{\it Claim 1. By Theorem 3.5 in (\cite{Os}, Claim 1), there are a
sequence $\{\Delta_i : i\in \mathbb{N}\}$ of pairwise disjoint
finite subsets of $X$, a sequence $\{\gamma_i: i\in \mathbb{N}\}$
of finite families of basis open sets in $\mathbb{R}^X$, and a
sequence $\{m_i: i\in \mathbb{N}\}\subseteq \mathbb{N}$ such that
the following conditions hold for every $i\in \mathbb{N}$:

$(a')$ $1\leq m_1$ and $m_i+2\leq m_{i+1}$;

$(b')$ $\overline{U}^{\mathbb{R}^X}\cap F_i=\emptyset$ and
$supp(U)\subseteq \bigcup_{j=1}^i \Delta_j$ for every $U\in
\gamma_i$;

$(c')$ if $f\in U\in \gamma_i$, then $|f(x)|\leq m_i$ for every
$x\in supp(U)$;

$(d')$ if $\varphi: A_i:=\bigcup\limits_{j=1}^i \Delta_j
\rightarrow \mathbb{R}$ is such that $\|\varphi\|_{A_i}\leq m_i$,
then there is a continuous function $f\in \bigcup \gamma_{i+1}$
such that $\|\varphi-f\|_{A_i}<\frac{1}{i}$}.

\medskip

Now we redefine the sequences in Claim 1 to make simpler and
clearer their usage in what follows.

By assumption the sequence $\{\Delta_n: n\in \mathbb{N}\}$
constructed in Claim 1 contains a strongly
$\mathrm{Coz}_{\delta}$-disjoint subsequence $\{\Delta_{n_k}: k\in
\mathbb{N}\}$, where $1<n_1<n_2<\dots$. Put

$R_1:=F_1$, $\Omega_1:=\bigcup\limits_{i=1}^{n_1-1}\Delta_i$,
$\mu_1:=\bigcup\limits_{i=1}^{n_1-1} \gamma_i$, $l_1:=m_{n_1-1}$,

and for every $k\in \mathbb{N}$, we define

$R_{2k}:=F_{n_k}$, $\Omega_{2k}:=\Delta_{n_k}$,
$\mu_{2k}:=\gamma_{n_k}$, $l_{2k}:=m_{n_k}$ and

$R_{2k+1}:=F_{n_k+1}$,
$\Omega_{2k+1}:=\bigcup\limits_{i=n_k+1}^{n_{k+1}-1} \Delta_{i}$,
$\mu_{2k+1}:=\bigcup\limits_{i=n_k+1}^{n_{k+1}-1}\gamma_{i}$, and

$l_{2k+1}:=m_{n_{k+1}-1}$.

It is clear that $\{R_i: i\in\mathbb{N}\}$ is an increasing
sequence of nowhere dense sets in $B^{st}_1(X)$ such that
$B^{st}_1(X)=\bigcup_{n\in \mathbb{N}} R_i$.

\medskip

{\it Claim 2. With respect to the sequence $\{R_i: i\in
\mathbb{N}\}$, the pairwise disjoint sequence $\{\Omega_i, i\in
\mathbb{N}\}$, and the sequences $\{\mu_i: i\in \mathbb{N}\}$ and
$\{l_i: i\in\mathbb{N}\}$ satisfy the following conditions $(i\in
\mathbb{N}):$

$(a)$ $1\leq l_1$ and $l_i+2\leq l_{i+1}$;

$(b)$ $\overline{U}^{\mathbb{R}^X}\cap R_i=\emptyset$ and
$supp(U)\subseteq \bigcup\limits_{j=1}^i \Omega_j$ for every $U\in
\mu_i$;

$(c)$ if $f\in U\in\mu_i$, then $|f(x)|\leq l_i$ for every $x\in
supp(U)$;

$(d)$ if $\varphi: A_i:=\bigcup_{j=1}^i \Omega_j\rightarrow
\mathbb{R}$ is such that $\|\varphi\|_{A_i}\leq l_i$, then there
is a continuous function $f\in \bigcup \mu_{i+1}$ such that
$\|\varphi-f\|_{A_i}<\frac{1}{i}$.

Moreover, the sequence $\{\Omega_{2i}: i\in \mathbb{N}\}$ is
strongly $\mathrm{Coz}_{\delta}$-disjoint.}

\medskip

{\it Proof of Claim 2.} By construction, $\{\Omega_i:i\in
\mathbb{N}\}$ is a sequence of pairwise disjoint finite subsets of
$X$, all families $\mu_i$ are finite, and the sequence
$\{\Omega_{2i}:i\in \mathbb{N}\}$ is strongly
$\mathrm{Coz}_{\delta}$-disjoint by the choice of
$\{\Delta_{n_k}:k\in \mathbb{N}\}$. The conditions (a) and (c) are
satisfied by $(a')$ and $(c')$, respectively. Since $F_i\subseteq
F_j$ for every $i\leq j$, we have $supp(U)\subseteq
\bigcup_{j=1}^i \Omega_j$ for every $U\in \mu_i$, and hence the
condition (b) holds true. The condition (d) is satisfied by the
definition of $\mu_i$ and the condition $(d')$ for the sets
$\gamma_i$. The claim is proved. \, \, \, \, \, \, \, \, \, $\Box$

\medskip

Since, by Claim 2, the sequence $\{\Omega_{2i}:i\in \mathbb{N}\}$
is strongly $\mathrm{Coz}_{\delta}$-disjoint, there is a disjoint
completely $\mathrm{Coz}_{\delta}$-additive family $\{W_x: x\in
\bigcup_{k\in \mathbb{N}} \Omega_{2k}\}$ of zero-sets of $X$ such
that $x\in W_x$ for all $x\in \bigcup_{k\in \mathbb{N}}
\Omega_{2k}$, and moreover, the sets

  $W(i):=\bigcup \{W_x: x\in \Omega_{2i}\}$ $(i\in \mathbb{N})$ are zero-sets in $X$.

\medskip

Let $F:=X\setminus \bigcup_{i\in \mathbb{N}} W(i)$. Then, $F\in
\mathrm{Zer}_{\sigma}(X)\cap \mathrm{Coz}_{\delta}(X)$.

\medskip

{\it Claim 3. By Theorem 3.5 in (\cite{Os}, Claim 4), there are a
sequence $\{f_i: i\in \mathbb{N}\}\subset B_1(X)$, a strictly
increasing sequence $\{r_i: i\in \mathbb{N}\}\subseteq \mathbb{N}$
with $r_1=1$, a sequence $\{b_i: i\in \mathbb{N}\}$, and a
sequence $\{U_i: i\in \mathbb{N}\}$ of basic open sets in
$\mathbb{R}^X$ such that

(e) $U_i\in \mu_{2r_i}$ for every $i\in \mathbb{N}$;

(f) $\|f_i\|_X\leq l_{2r_i}$ for every $i\in \mathbb{N}$;

(g) $\|f_{i+1}-f_i\|_{C_i}\leq \frac{1}{2^{b_i}}$ for every $i\in
\mathbb{N}$, where $C_i:=F\cup \bigcup_{j=1}^{r_i} W(j)$;

(h) $f_j\in U_i$ for every $j\geq i\geq 1$;

(k) $b_{i+1}-b_{i}\rightarrow \infty$.}

\medskip

It follows from Claim 4(g) that the sequence $\{f_i: i\in
\mathbb{N}\}\subseteq B_1(X)$ converges pointwise to some function
$f: X\rightarrow \mathbb{R}$.


Since $F\in \mathrm{Zer}_{\sigma}(X)$, then $F=\bigcup_{i\in
\mathbb{N}} Q_i$ where $Q_i\subseteq Q_{i+1}$ and $Q_i$ is a
zero-set of $X$ for every $i\in \mathbb{N}$.

Since $\{Q_i,\{W_z: z\in \bigcup\limits_{j=1}^i \Omega_{2i}\}\}$
is a pairwise disjoint finite family of zero-sets of $X$,  one can
define a continuous function
 $g_{i}: X\rightarrow
\mathbb{R}$ such that $g_{i}(x)=0$ for $x\in Q_i$ and
$g_{i}(x)=f(z)$ for every $x\in W_z$ where $z\in
\bigcup\limits_{j=1}^i \Omega_{2i}$.

The sequence $\{g_i: i\in \mathbb{N}\}\subseteq C_p(X)$ converges
pointwise to the function $g: X\rightarrow \mathbb{R}$ where

$$ g(x):=\left\{
\begin{array}{lcl}
0, \, \, \, \, \, if \, \, \, \, x\in F;\\
f(z), \, \, \, \, \, if \, \, \, \, x\in W_z \, \,  for \, some \, \, z\in \bigcup\limits_{i\in \mathbb{N}} \Omega_{2i}. \\
\end{array}
\right.
$$

Note that $g\in B^{st}_1(X)$.

By Claim 3(h),  $g\in \overline{U_i}$ for each $i\in \mathbb{N}$.
But, then, by Claim 2(b), $g\not\in \bigcup\limits_{i=1}^{\infty}
F_i=B^{st}_1(X)$, a contradiction.

\end{proof}

\medskip

Let  $\mathfrak{b}=\min\{|X|:$ $X$ has a countable pseudocharacter
but $X$ is not a $\lambda$-space$\}$ (see \protect\cite{handbook},
p.149).

\medskip

By Theorem \ref{th23}, Proposition 3.15 and Example 4.8 in
\cite{BG}, we get the following results:

$\bullet$ If $X$ is metrizable and $|X|<\mathfrak{b}$ then
$B^{st}_1(X)$ is Baire.

$\bullet$ (CH) If $X$ is a Luzin set then $B^{st}_1(X)$ is meager.

$\bullet$ It is consistent with $ZFC$ there is a zero-dimensional
metrizable separable space $X$ with $|X|=\mathfrak{b}$  such that
$B_1^{st}(X)$ is Baire but $B_1^{st}(X)$ is not Choquet (see the
definition of a {\it Choquet} space in \cite{BG}).

\medskip


\bibliographystyle{model1a-num-names}
\bibliography{<your-bib-database>}

\end{document}